\documentclass[11pt,a4paper]{article}
\usepackage{color}
\usepackage[colorlinks=true]{hyperref}
\usepackage[margin=1in]{geometry}
\usepackage{fullpage,mathrsfs,graphicx,framed,color,amssymb,amsmath,amsthm}
\usepackage[boxed, noline, ruled, linesnumbered]{algorithm2e}
\usepackage{epsfig, graphicx}
\usepackage{latexsym,amsfonts,amsbsy,amssymb}
\usepackage{amsmath,amsthm}
\usepackage{enumerate}
\usepackage{cleveref}
\usepackage{changes} 
\makeatletter 

\@addtoreset{equation}{section}
\makeatother 
\allowdisplaybreaks 
\newtheorem{theorem}{Theorem}[section]

\allowdisplaybreaks 

\numberwithin{equation}{section}

\begin{document}
\title{Tensor Neural Network Interpolation and Its Applications\footnote{This work 
was supported by the National Key Research and Development Program 
of China (2023YFB3309104), National Natural Science Foundations of 
China (NSFC 1233000214),  National Key Laboratory of Computational Physics 
(No. 6142A05230501),  Beijing Natural Science Foundation (Z200003), 
National Center for Mathematics and Interdisciplinary Science, CAS.}}
\author{Yongxin Li\footnote{School of Statistics and Mathematics, 
Central University of Finance and Economics, 
No.39, Xueyuan Nanlu, Beijing, 100081, China (2018110075@email.cufe.edu.cn)},\ \ \ 
Zhongshuo Lin\footnote{LSEC,
Academy of Mathematics and Systems Science,
Chinese Academy of
Sciences, No.55, Zhongguancun Donglu, Beijing 100190, China, and School of
Mathematical Sciences, University of Chinese Academy
of Sciences, Beijing, 100049 (linzhongshuo@lsec.cc.ac.cn)},\ \ \
Yifan Wang\footnote{LSEC,
Academy of Mathematics and Systems Science,
Chinese Academy of
Sciences, No.55, Zhongguancun Donglu, Beijing 100190, China, and School of
Mathematical Sciences, University of Chinese Academy
of Sciences, Beijing, 100049 (wangyifan@lsec.cc.ac.cn)}\ \ \
and \ \ Hehu Xie\footnote{LSEC,
Academy of Mathematics and Systems Science,
Chinese Academy of
Sciences, No.55, Zhongguancun Donglu, Beijing 100190, China, and School of
Mathematical Sciences, University of Chinese Academy
of Sciences, Beijing, 100049 (hhxie@lsec.cc.ac.cn)}}

\date{}
\maketitle

\begin{abstract}
Based on tensor neural network, we propose an interpolation method for 
high dimensional non-tensor-product-type functions. This interpolation scheme 
is designed by using the tensor neural network based machine learning method. 
This means that we use a tensor neural network to approximate high dimensional 
functions which has no tensor product structure. In some sense, 
the non-tenor-product-type high dimensional function is transformed to the 
tensor neural network which has tensor product structure. 
It is well known that the tensor product structure can bring the possibility to 
design highly accurate and efficient numerical methods for dealing with high 
dimensional functions. In this paper, we will concentrate on computing 
the high dimensional integrations and solving high dimensional partial differential 
equations. The corresponding numerical methods and numerical examples  will be provided 
to validate the proposed tensor neural network interpolation.


\end{abstract}

\section{Introduction}
There exist many high-dimensional problems such as quantum mechanics, 
statistical mechanics, and financial engineering in modern sciences and engineering. 
Traditional numerical methods like finite difference, finite element, 
and spectral methods are typically confined to solving low-dimensional problems. 
However, applying these classical methods to high-dimensional problems 
will encounter the so-called ``curse of dimensionality''. 

Due to its universal approximation property, the fully connected neural network (FNN) 
is the most widely used architecture to build the functions for solving 
high-dimensional PDEs. There are several types of well-known FNN-based methods 
such as deep Ritz \cite{EYu}, deep Galerkin method \cite{DGM}, 
PINN \cite{RaissiPerdikarisKarniadakis}, and weak adversarial networks \cite{WAN}
for solving high-dimensional PDEs by designing different loss functions. 
Among these methods, the loss functions always include computing 
high-dimensional integrations for the functions defined by FNN. 
For example, the loss functions of the deep Ritz method, deep Galerkin method 
and weak adversarial networks method require computing the integrations 
on the high-dimensional domain for the functions constructed by FNN. 
Direct numerical integration for the high-dimensional functions also 
meets the ``curse of dimensionality''. Always, the high-dimensional integration 
is computed using the Monte-Carlo method along with some sampling tricks \cite{EYu,HanZhangE}. 
Due to the low convergence rate of the Monte-Carlo method, the solutions obtained 
by the FNN-based numerical methods are challenging to achieve high accuracy and 
stable convergence process. This means that the Monte-Carlo method decreases 
computational work in each forward propagation while decreasing the simulation 
accuracy, efficiency and stability of the FNN-based numerical methods 
for solving high-dimensional PDEs. 

Recently, we propose a type of tensor neural network (TNN) and the corresponding
machine learning method to solve high-dimensional problems with high accuracy.
The most important property of TNN is that the corresponding high-dimensional functions
can be easily integrated with high accuracy and high efficiency. Then the deduced machine
learning method can arrive high accuracy for solving high-dimensional problems.
The reason is that the high dimensional integration of TNN in the loss functions 
can be transformed into one-dimensional integrations which can be computed 
by the classical quadrature schemes with high accuracy.
The TNN based machine learning method has already been used to solve high-dimensional
eigenvalue problems and boundary value problems based on the Ritz type of loss
functions \cite{WangJinXie}. Furthermore,  in \cite{WangXie}, the multi-eigenpairs 
can also been computed with machine learning method which is designed by combining the TNN 
and Rayleigh-Ritz process. Furthermore, with the help of TNN, the 
a posteriori error estimator can be adopted as the loss function of the machine learning method 
for solving high dimensional boundary value problems and eigenvalue problems \cite{WangLinLiaoLiuXie}.   
So far, the TNN based machine learning method has shown good ability for solving 
high dimensional problems which having only tensor product coefficients and source terms. 
The aim and main contribution of this paper is to design a type of TNN based interpolation method 
for the non-tensor-product-type functions. Then using this interpolation method to approximate 
the non-tensor-product coefficients and then modify the non-tensor-product-type problems to the 
corresponding tensor product one. Then the TNN based machine learning method can be 
adopted to solve the deduced tensor product problems with high accuracy.

An outline of the paper goes as follows. 
In Section \ref{Section_Integration_Accuracy}, 
we investigate the dependence of the accuracy for machine learning methods 
on the accuracy of integration through a numerical experiment. 
Section \ref{Section_TNN} is devoted to introducing the TNN structure and the corresponding 
approximation property. In Section \ref{Section_TNN_Interpolation}, 
the TNN based interpolation will be proposed to approximate high-dimensional functions.
In Section \ref{Section_Integration}, we introduce the way to 
compute the high dimensional integrations of non-tensor-product-type
functions by using TNN interpolation. Section \ref{Section_TNN_PDEs} is devoted 
to introducing the application of TNN interpolation to solving high dimensional 
partial differential equations with non-tensor-product-type coefficients and source terms. 
Section \ref{Section_Numerical_Examples} gives some numerical examples to validate
the accuracy and efficiency of the proposed TNN based machine learning method.
Some concluding remarks are given in the last section.
\section{Integration accuracy and machine learning accuracy}\label{Section_Integration_Accuracy} 
As mentioned above, the loss functions in common machine 
learning methods for solving PDEs always include high dimensional 
integrations for the functions defined by neural networks. 
In this section, we investigate the dependence of the accuracy of 
machine learning method on the integration accuracy. 
The aim here is to explore the importance of the integration 
accuracy for the machine learning method.  

We do the numerical investigations of the machine learning accuracy 
for the Deep Ritz and PINN methods with different integration schemes. 
To be specific, we consider the following eigenvalue problems: 
Find $(\lambda,u)\in \mathbb R \times H_0^1(\Omega)$ such that 
\begin{eqnarray*}
\left\{
\begin{aligned}
-\Delta u&=\lambda u,\ \ \ &x\in& \Omega,\\
u&=0,\ \ \ &x\in&\partial \Omega,
\end{aligned}
\right.
\end{eqnarray*}
where $\Omega := [0,1]^2$. We consider this two dimensional problem so that computing the integrations 
included in the loss function using classical Gauss quadrature scheme becomes feasible. 
Let $\phi(x;\theta)$ denote the neural network. 
Here, we check and compare the numerical performances of the following four strategies:
\begin{enumerate}
\item (RitzMC): Use the Deep Ritz method with the following loss function
\begin{eqnarray}\label{loss_DRM}
\mathcal L_{\rm DRM}(x;\theta) 
= \frac{\int_\Omega|\nabla \phi(x;\theta)|^2 dx}{\int_\Omega |\phi(x;\theta)|^2 dx}
\end{eqnarray}
and compute the included two dimensional integrations with Monte Carlo integration method.
\item (RitzGauss): Use the Deep Ritz method with loss function (\ref{loss_DRM}) 
and compute the included two dimensional integrations with Gauss quadrature method.
\item (PINNMC): Use the following loss function based on the idea of PINN 
method and Deep Galerkin method
\begin{eqnarray}\label{loss_PINN}
\mathcal L_{\rm PINN}(x;\theta)&=&\|\Delta \phi(x;\theta)
+\lambda(x;\theta)\phi(x;\theta)\|_{L^2(\Omega)}^2 \nonumber\\
&=&\left\|\Delta \phi(x;\theta)+\frac{\int_{\Omega}|\nabla\phi(x;\theta)|^2dx}{\int_\Omega |\phi(x;\theta)|^2dx}\phi(x;\theta)\right\|_{L^2(\Omega)}^2.
\end{eqnarray}
and compute the two dimensional integrations with Monte Carlo integration method.
\item (PINNGauss): Use the loss function (\ref{loss_PINN}) and compute 
the two dimensional integration with Gauss quadrature method.
\end{enumerate}
To guarantee the fairness of the numerical comparison, 
we use the same neural network structures in all four experiments. 
The fully connected neural network with three hidden layers, with each layer having 50 neurons, 
is adopted here. 
The sine function is chosen as the activation function. 
In order to guarantee the homogeneous boundary condition, 
we multiply $\phi(x;\theta)$ by $x_1(1-x_1)x_2(1-x_2)$, 
and still denote the resulting function as $\phi(x;\theta)$ 
for notational convenience. The same number of quadrature points is used 
to compute the two dimensional integrations included in the loss functions. 
Specifically, for Gauss quadrature method, we decompose the 
interval $[0,1]$ into $20$ subintervals and choose $8$ Gauss quadrature points 
in each subinterval to obtain the composite quadrature points 
set $X$ of one dimension. Then, we generate the two dimensional 
quadrature points using the Cartesian product of $X$ and itself, 
i.e., $X \times X$, which contains $25,600$ two dimensional quadrature 
points in total. As for the Monte Carlo integration method, 
we sample 25,600 uniformly distributed points in $[0,1]^2$ 
in each iteration step. The Adam optimizer with learning rate 0.0003 
is adopted for 500,000 iteration steps in all four experiments.

\begin{figure}[htb!]
\centering
\includegraphics[width=9cm,height=9cm]{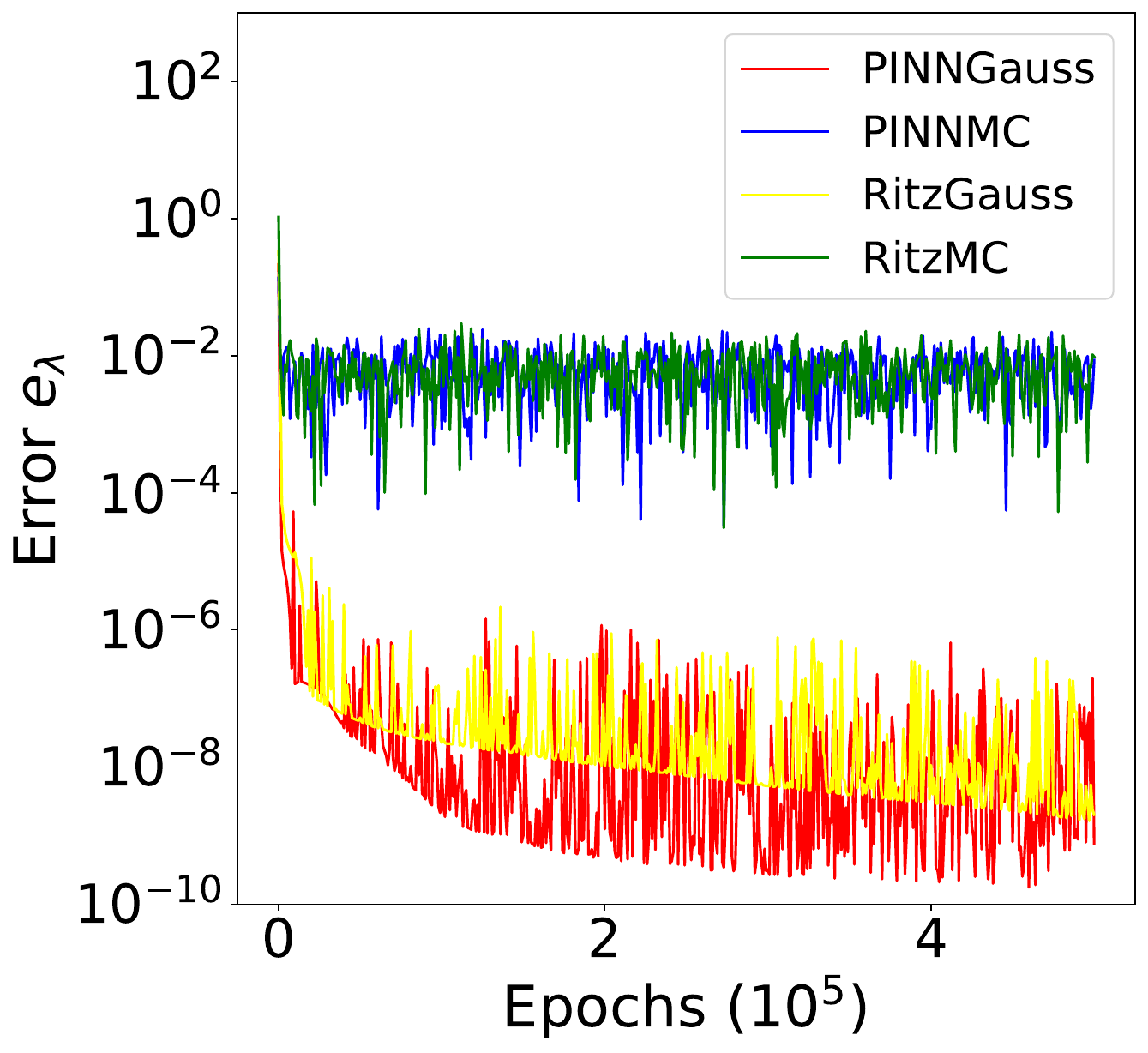}
\caption{Relative errors of four methods during the training 
process.}\label{fnn_laplace_errorE}
\end{figure}

We record the network model, and compute the relative error
\[
e_{\lambda}:= \frac{|\lambda - \lambda^*|}{\lambda^*}    
\]
of the corresponding approximate eigenvalue $\lambda$ for the exact 
eigenvalue $\lambda^*$ after each $100$ steps. Figure \ref{fnn_laplace_errorE} 
shows the relative errors of four methods during the training process.
It can be easily observed that the performances of two different 
loss functions (\ref{loss_DRM}) and (\ref{loss_PINN}) are similar, 
while what really matters is the choice of the numerical integration method. 
To be specific, with the high-accuracy Gauss quadrature method, 
the RitzGauss and PINNGauss methods are able to achieve obviously higher 
accuracy than that of the RitzMC and PINNMC methods. 

The results above give a hint that the accuracy of computing the integrations 
included in the loss functions indeed plays a significant role for the accuracy of 
the machine learning methods for solving PDEs. 
And therefore, it's necessary for the machine learning methods 
to guarantee the computation accuracy, 
such as the integration accuracy during the iteration process. 
Notice that we consider the two dimensional problem in this section 
since the two dimensional Gauss quadrature method can be leveraged in this case. 
However, in high dimensional cases, using Gauss quadrature method 
to guarantee the integration accuracy is not realistic 
and feasible due to the curse of dimensionality. 
And this motivates us to design the tensor neural network structure 
and TNN based machine learning methods 
\cite{WangJinXie,WangLinLiaoLiuXie,WangLiaoXie}, 
which assure the high dimensional integration accuracy within acceptable 
computational complexity, and therefore achieve high accuracy 
in solving high dimensional PDEs.

\section{Tensor neural network architecture}\label{Section_TNN} 
\begin{figure}[htb!]
\centering
\includegraphics[width=12cm,height=9cm]{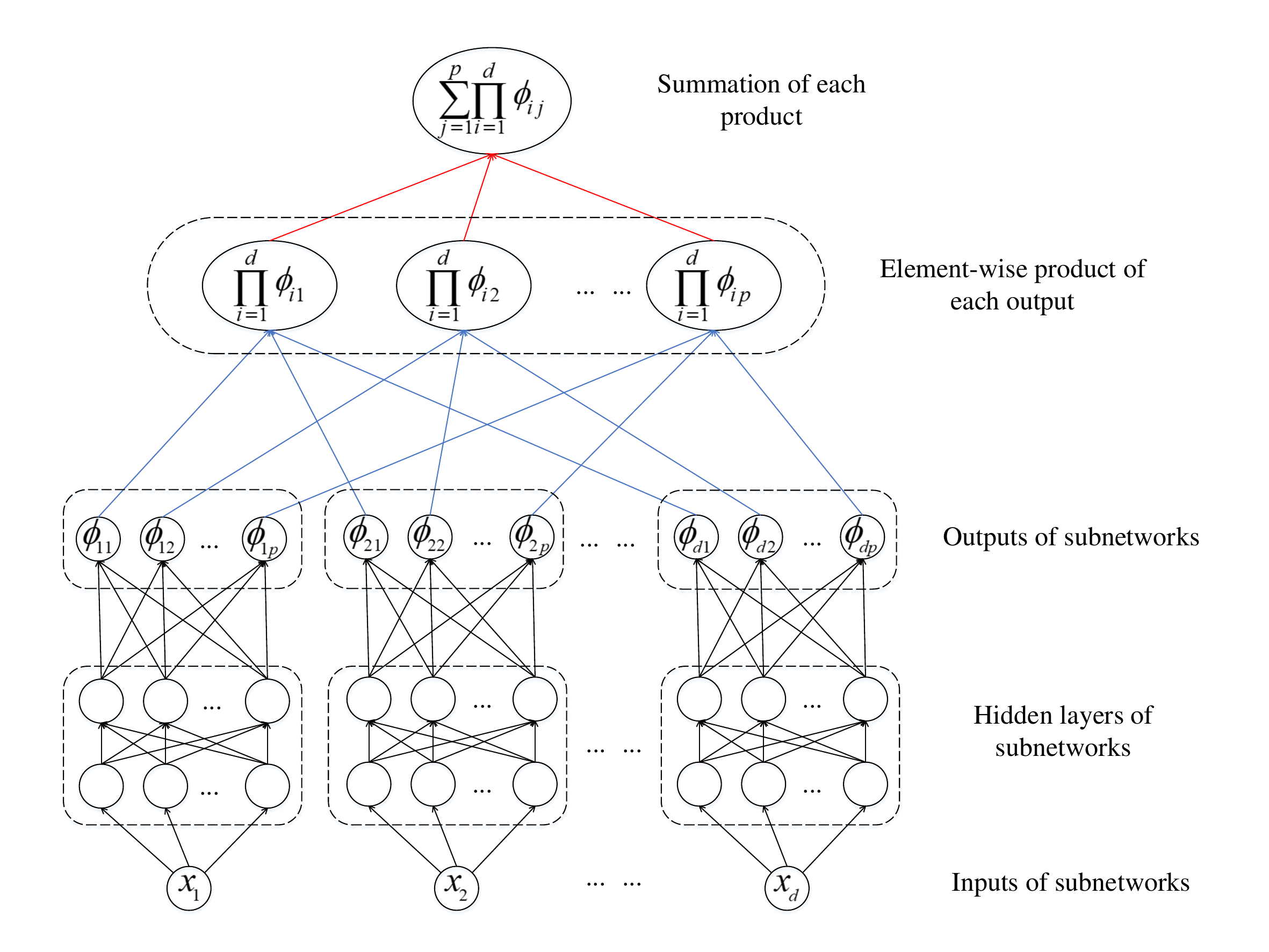}
\caption{Architecture of TNN. Black arrows mean linear transformation
(or affine transformation). Each ending node of blue arrows is obtained by taking the
scalar multiplication of all starting nodes of blue arrows that end in this ending node.
The final output of TNN is derived from the summation of all 
starting nodes of red arrows.}\label{TNNstructure}
\end{figure}

In this section, we introduce the architecture of TNN and 
its approximation property which have been discussed and 
investigated in \cite{WangJinXie}. The TNN is a neural network 
of low-rank structure, which is built by the tensor product of 
several one-dimensional input and multi-dimensional output 
subnetworks. Due to the low-rank structure of TNN, an efficient 
and accurate quadrature scheme can be designed for the
TNN-related high-dimensional integrations such as the inner 
product of two TNNs. In \cite{WangJinXie}, we introduce TNN 
in detail and propose its numerical integration scheme with
the polynomial scale computational complexity of the dimension.
For each $i=1,2,\cdots,d$, we use $\Phi_i(x_i;\theta_i)
=(\phi_{i,1}(x_i;\theta_i),\phi_{i,2}(x_i;\theta_i), 
\cdots,\phi_{i,p}(x_i;\theta_i))$ to denote a subnetwork that 
maps a set $\Omega_i\subset\mathbb R$ to $\mathbb R^p$,
where $\Omega_i,i=1,\cdots,d,$ can be a bounded interval $(a_i,b_i)$, 
the whole line $(-\infty,+\infty)$ or the half line $(a_i,+\infty)$.
The number of layers and neurons in each layer, the selections of 
activation functions and other
hyperparameters can be different in different subnetworks.
In this paper, in order to improve the numerical stability further, 
the TNN is defined as follows:
\begin{eqnarray}\label{def_TNN}
\Psi(x;\Theta)=\sum_{j=1}^pc_j\widehat\phi_{1,j}(x_1;\theta_1)
\widehat\phi_{2,j}(x_2;\theta_2)
\cdots\widehat\phi_{d,j}(x_d;\theta_d)
=\sum_{j=1}^pc_j\prod_{i=1}^d\widehat\phi_{i,j}(x_i;\theta_i),
\end{eqnarray}
where $c=\{c_j\}_{j=1}^{p}$ is a set of trainable parameters, 
$\Theta=\{c,\theta_1,\cdots,\theta_d\}$
denotes all parameters of the whole architecture.
For $i=1,\cdots,d,j=1,\cdots,p$, $\widehat\phi_{i,j}(x_i,\theta_i)$ 
is a normalized functions as follows:
\begin{eqnarray}\label{eq_phi_normed}
\widehat\phi_{i,j}(x_i,\theta_i)
=\frac{\phi_{i,j}(x_i,\theta_i)}{\|\phi_{i,j}(x_i,\theta_i)\|_{L^2(\Omega_i)}}.
\end{eqnarray}
The TNN architecture (\ref{def_TNN}) and the one defined in \cite{WangJinXie} 
are mathematically equivalent, but (\ref{def_TNN}) has better numerical 
stability during the training process.
Figure \ref{TNNstructure} shows the corresponding architecture of TNN.
From Figure \ref{TNNstructure} and numerical tests, we can find the parameters 
for each rank of TNN are correlated by the FNN, which guarantee the 
stability of the TNN-based machine learning methods. This is also an important 
difference from the tensor finite element methods.

In \cite{WangJinXie} we introduce and prove the following approximation result 
to the functions in the space $H^m(\Omega_1\times\cdots\times\Omega_d)$ 
under the sense of $H^m$-norm.

\begin{theorem}\label{theorem_approximation}
Assume that each $\Omega_i$ is an interval in $\mathbb R$ for $i=1, \cdots, d$, $\Omega=\Omega_1\times\cdots\times\Omega_d$,
and the function $f(x)\in H^m(\Omega)$. Then for any tolerance $\varepsilon>0$, there exist a
positive integer $p$ and the corresponding TNN defined by (\ref{def_TNN})
such that the following approximation property holds
\begin{equation}\label{eq:L2_app}
\|f(x)-\Psi(x;\theta)\|_{H^m(\Omega)}<\varepsilon.
\end{equation}
\end{theorem}
The TNN seems rather simple, but it actually has a surprising rich structure.
The motivation for employing the TNN architectures is to provide 
high accuracy and high efficiency in calculating variational forms 
of high-dimensional problems in which high-dimensional
integrations are included. The TNN itself can approximate functions in 
Sobolev space with respect to $H^m$-norm. For more information about the rank 
estimates, please refer to \cite{WangJinXie}. 
The major contribution of this paper is to propose a TNN-based interpolation method, 
to approximate high dimensional functions. Its applications for computing 
high dimensional integrations and solving high dimensional 
partial differential equations also shows the advantages of TNN.

\section{Tensor neural network interpolation}\label{Section_TNN_Interpolation}
In this section, we introduce the tensor neural network interpolation method 
for the desired function $f(x)$. In this article, we say that a function $g(x)$ 
is of tensor product type, or is a tensor-product-type function, if $g(x)$ 
can be expressed in the form of
\begin{equation}\label{T-P-type}
g(x) = \sum_{j=1}^p g_{1,j}(x_1)g_{2,j}(x_2)\cdots g_{d,j}(x_d) = \sum_{j=1}^{p} 
\prod_{i=1}^{d}g_{i,j}(x_i)
\end{equation}
for some one-dimensional functions $g_{i,j}(x_i)$. Recall that TNN function 
$\Psi(x;\Theta)$ defined by (\ref{def_TNN}) is of tensor product type with 
$\widehat\phi_{i,j}(x_i;\theta_i)$ being one-dimensional-input neural network. 
In the previous section, we demonstrate the importance of integration accuracy 
on the accuracy of machine learning method for solving PDEs. 
When the dimension of the problem is high, computing the high dimensional 
integration within a satisfactory accuracy under limited computing resources 
is difficult, unless the integrated function possesses some desirable properties.  
In \cite{WangJinXie}, the efficient and accurate integration scheme for 
a function of tensor product type is designed. The core idea is to transform 
the high dimensional integration into the product of one dimensional one, 
and use classical numerical integration scheme such as Gauss quadrature 
scheme to compute these one dimensional integration. However, in practice, 
we also encounter some scenarios where the integrated function is not of 
tensor product type. For example, when we use machine learning method to 
solve some high dimensional problems with non-tensor-product-type coefficients 
or source terms, the computation of the loss function is always tricky.

The aim of the interpolation is to find 
a tensor neural network $\Psi(x,\theta)$ to approximate the integrated 
function $f(x)$ (possibly non-tensor-product-type)
in some sense so that we can use the approximated function to compute 
the high dimensional integration of $f(x)$. Note that finding an approximation for 
high dimensional functions is generally not easier than computing the integrations 
of these high dimensional functions. However, in machine learning method for PDEs, 
there requires lots of high dimensional integrations computation.
Then finding an approximated function that can be integrated 
with high efficiency and accuracy can improve the accuracy of 
machine learning method and speed up the learning process, 
since we only do the TNN interpolation once and use the approximate 
function to replace the original one during the whole training process. 
Therefore, the interpolation trick is of great importance for TNN 
based machine learning methods to be applied again 
in the case of non-tensor-product-type coefficients 
or source terms, due to the tensor product structure of the 
TNN approximation function.

It is a common idea to use an easy-to-integrated function to approximate 
the integrated function. For example, in one dimensional integration, 
classical quadrature schemes are designed 
by using polynomial interpolation for the integrated functions. 
Especially, the idea behind the well known Gauss quadrature schemes  
are based on the special choice of the interpolation points. 
In this paper, we use TNN to approximate the integrated function 
due to its universal approximation ability and the fact that TNN function 
can be integrated in an efficient and accurate way \cite{WangJinXie}.

The machine learning method is adopted to finish this approximating task. 
To be specific, at each iteration $\ell$, we obtain a bunch of training points 
$x_k^{(\ell)} := (x_{k,1}^{(\ell)}, \cdots, x_{k,d}^{(\ell)})^\top$, 
$k=1, \cdots, K$ according to some sampling rules, and minimize the squared loss
\begin{equation}
\begin{aligned}
L_\ell(\Theta) :=& \sum_{k=1}^K \left(\Psi(x_k^{(\ell)},\Theta)-f(x_k^{(\ell)})\right)^2\\
=& \sum_{k=1}^K \left(\sum_{j=1}^pc_j\prod_{i=1}^d
\widehat\phi_{i,j}(x_{k,i}^{(\ell)};\theta_i)-f(x_k^{(\ell)})\right)^2,
\end{aligned}
\end{equation}
to obtain the desired network parameters $\Theta=\{c,\theta_1,\cdots,\theta_d\}$. 
This procedure is repeated $M$ times until we obtain good enough results on 
the validation scheme, such as the accuracy on the test points set.

In the optimization process, we split the parameters into 
two groups $\{c\}$ and $\{\theta_1,\cdots,\theta_d\}$. 
The parameter $c$ can be regarded as the linear coefficients on 
the $p$-dimensional subspace $V_p^{(t)} :=\text{span}\left\{\varphi_{j}(x;\theta^{(t)})
:=\prod_{i=1}^{d}\widehat\phi_{i,j}\left(x_{i};\theta_{i}^{(t)}\right)\right\}$. 
And therefore, we only need to solve a linear equation to obtain the optimal 
coefficient $c$ on the current subspace $V_p^{(t)}$ in the sense of the 
squared loss on training points. 
Using the optimal coefficient $c$, we update the network parameters $\{\theta_j\}$ 
by minimizing the loss function with some optimization algorithm. 
The TNN interpolation method to obtain an approximation 
for a given target function is defined in Algorithm \ref{Algorithm_Interpolation}. 
\begin{algorithm}[htb!]
\caption{TNN interpolation method}\label{Algorithm_Interpolation}
\SetKwInOut{HP}{Hyper-parameters}
\KwIn{Target function $f(x)$, TNN function $\Psi(x;\Theta)$ defined by (\ref{def_TNN}), 
initial model parameters $\Theta$, domain $\Omega$.}
\KwOut{Learned approximate TNN function $\Psi(x;\Theta^*)$}
\HP{Number of total training iterations $M$, number of training points in each iteration $K$, 
number of optimization steps $N$ for each training points set, hyperparameters 
of optimization algorithm such as step size $\gamma$.}
\BlankLine
\For{$\ell \leftarrow 1$ \KwTo $M$}{
    Sample $x_k^{(\ell)} \in \Omega$, $k=1,\dots,K$ according to some sampling rules \\
    \For{$t \leftarrow 1$ \KwTo $N$}{
        Assemble matrix $A^{(t)}$ and vector $B^{(t)}$, where the entries are defined as follows
        \begin{equation*}
            \begin{aligned}
                A_{m,n}^{(t)}&= \sum_{k=1}^{K} \prod_{i=1}^{d} \widehat\phi_{i,m}(x_{k,i}^{(\ell)};\theta_i^{(t-1)}) \prod_{i=1}^{d} \widehat\phi_{i,n}(x_{k,i}^{(\ell)};\theta_i^{(t-1)}), \ 1\leq m,n\leq p, \\
                B_{m}^{(t)} &= \sum_{k=1}^{K} f(x_k^{(\ell)}) \prod_{i=1}^{d} \widehat\phi_{i,m}(x_{k,i}^{(\ell)};\theta_i^{(t-1)}), \ 1\leq m\leq p.
            \end{aligned}
        \end{equation*}

        Solve the following linear equation to obtain the solution $c \in \mathbb{R}^p$:
        \begin{equation*}
            A^{(t)}c = B^{(t)},
        \end{equation*}
        and update the coefficient parameter as $c^{(t)} = c$.
        
        Compute the loss function
\begin{equation*}
\begin{aligned}
\mathcal L_\ell^{(t)}(\theta^{(t-1)}) = \sum_{k=1}^K \left(\sum_{j=1}^pc_j^{(t)}\prod_{i=1}^d\widehat\phi_{i,j}\left(x_{k,i}^{(\ell)};\theta_i^{(t-1)}\right)
-f(x_k^{(\ell)})\right)^2.
\end{aligned}
\end{equation*}

        Use an optimization step to update the neural network parameters of TNN as follows:
        \begin{equation*}
            \theta^{(t)} = \theta^{(t-1)} - \gamma\frac{\partial\mathcal L_\ell^{(t)}}{\partial\theta}(c^{(t)},\theta^{(t-1)}).
        \end{equation*}
    }
}
\end{algorithm}
For different target function $f(x)$ and computing domain $\Omega$,  
we can literally design all the details as we need, 
including the parameters, the validation scheme and so on. 
Recall that our primary motivation is to obtain an optimal function, 
which can be integrated with high efficiency and accuracy,  
to replace the non-tensor-product coefficient in the PDEs 
and then use the TNN based machine learning method to solve the concerned PDEs. 
Therefore all the interpolation work can be done off-line, 
and we can repeat the procedure until we obtain an satisfactory 
appropriate function as we need.

\section{TNN interpolation for high dimensional integration}\label{Section_Integration}
In this section, based on the TNN interpolation, we design a type of machine learning method 
for computing high dimensional integrations. As an example, we are concerned with 
the following high dimensional integration for the $f(x)$ on $\Omega\subset \mathbb R^d$
\begin{eqnarray}\label{Integration}
I = \int_\Omega f(x)dx.
\end{eqnarray}
Actually, the process is direct and simple. First, we produce the TNN interpolation 
$\Psi(x,\theta)$ to approximate the target function $f(x)$ 
based on the machine learning process defined by Algorithm \ref{Algorithm_Interpolation}. 
Then, with the help of TNN interpolation  $\Psi(x,\theta)$, 
we can compute the following integration as an approximation 
to the high dimensional integration (\ref{Integration})
\begin{eqnarray}
I \approx \int_\Omega \Psi(x,\theta)dx. 
\end{eqnarray}

The method here can be extended to compute the numerical integrations of polynomial 
composite functions of TNN and their derivatives. 
For more information about computing high dimensional integrations of TNN functions, 
please refer to \cite{WangJinXie}. For example, the method here can be used to compute the following 
integrations
\begin{eqnarray}
\int_\Omega |f|^2dx,\ \ \ \int_{\Omega}|\nabla f|^2dx,\ \ \ 
\int_{\Omega}|-\Delta f|^2dx.
\end{eqnarray}
About the way and complexity for computing numerical integrations of 
polynomial composite functions of TNN and their derivatives, 
please refer to  \cite{WangJinXie}.  
For example, we have the following way to compute $\int_\Omega \Psi dx$ and  $\int_\Omega |\Psi|^2dx$
\begin{eqnarray*}
\int_\Omega \Psi(x, \theta)dx &=& \sum_{j=1}^p\prod_{i=1}^d\int_{\Omega_i} \phi_{i, j}(x_i, \theta_i)dx_i
\approx\sum_{j=1}^p\prod_{i=1}^d
\left(\sum_{n_i=1}^{N_i}w_i^{(n_i)}\phi_{i,j}(x_i^{(n_i)})\right),\\
\int_\Omega|\Psi|^2dx&=&\sum_{j=1}^p\sum_{k=1}^p\prod_{i=1}^d
\left(\int_{\Omega_i}\phi_{i,j}(x_i)\phi_{i,k}(x_i)dx_i\right)\nonumber\\
&\approx&\sum_{j=1}^p\sum_{k=1}^p\prod_{i=1}^d
\left(\sum_{n_i=1}^{N_i}w_i^{(n_i)}\phi_{i,j}(x_i^{(n_i)})
\phi_{i,k}(x_i^{(n_i)})\right),
\end{eqnarray*}
where $\Omega=\Omega_1\times\Omega_d$, $(x_i^{(n_i)},w_i^{(n_i)})$ denotes the 
quadrature points and weights on the domain $\Omega_i$, $n_i=1, \cdots, N_i$.

The method here give a new view to understand the numerical interpolation 
and numerical integration. Different from Monte-Carlo based schemes, 
it give another way to compute the high dimensional numerical 
integrations. More about this topic, please refer to \cite{FengZhong, HuaWang,KuoSchwabSloan,XuZhou}. 

\section{TNN interpolation for solving high dimensional PDEs}\label{Section_TNN_PDEs}
In this section, we introduce the application of TNN interpolation method for solving 
high dimensional elliptic problem with non-tensor-product-type coefficients and source term. 
The numerical scheme here is designed based on the combination of 
the TNN interpolation method and TNN-based machine learning method \cite{WangLinLiaoLiuXie}.
We assume the physical domain $D=D_1\times\cdots\times D_d$ with $D_i=[a_i, b_i]$
for $i=1, \cdots, d$. It is noted that the tensor product structure plays
an important role in reducing the dependence on dimensions in numerical integration.
We will find that the high accuracy of the high-dimensional integrations of
TNN functions leads to the corresponding machine learning method has high accuracy
for solving the high-dimensional parametric elliptic equations. 

As an example, we consider the following second order elliptic problem: 
Find $\bar u\in H_0^1(\Omega)$ such that 
\begin{eqnarray}\label{Laplace_Problem}
\left\{
\begin{array}{rcl}
-{\rm div}(a(x)\nabla \bar u) +b(x)\bar u&=&f(x),\ \ \ \ {\rm in}\ \Omega,\\
\bar u&=&0, \ \ \quad \quad {\rm on}\ \partial\Omega,
\end{array}
\right.
\end{eqnarray}
where the coefficients $a(x)\geq a_{\min}>0$ and $b(x)\geq 0$, 
the source term $f\in L^2(\Omega)$ may be non-tensor-product-type functions.  

First, we use the TNN interpolation method defined in Algorithm \ref{Algorithm_Interpolation} 
to get the tensor neural network functions $\widehat a(x)$, $\widehat b(x)$ 
and $\widehat f(x)$, respectively. Then the concerned equation  (\ref{Laplace_Problem}) 
can be modified to the following approximate one: Find $u\in H_0^1(\Omega)$ such that 
\begin{eqnarray}\label{Laplace_Problem_Modified}
\left\{
\begin{array}{rcl}
-{\rm div}(\widehat{a}(x)\nabla u) + \widehat b(x)u&=&\widehat f(x),\ \ \ \ {\rm in}\ \Omega,\\
u&=&0, \ \ \ \ {\rm on}\ \partial\Omega.
\end{array}
\right.
\end{eqnarray}
The solution $u$ here can be regarded as the approximation to the exact solution $\bar u$ 
of the equation (\ref{Laplace_Problem}). 

For designing the TNN based machine learning method, we build the following TNN function
\begin{eqnarray}\label{tab:TNN_train}
\Psi(x)=\sum\limits_{j=1}^p c_j\prod_{i=1}^d\widehat\phi_{i,j}(x_i),
\end{eqnarray}
where $x=[x_1, \cdots, x_d]^\top$. 
In order to deal with the boundary condition, following \cite{WangJinXie},
for $i=1,\cdots, d$, the $i$-th subnetwork $\psi_i(x_i;\theta_i)$ is defined as follows:
\begin{eqnarray*}\label{bd_subnetwork}
\phi_i(x_i)&:=&(x_i-a_i)(b_i-x_i)\widehat\phi_i(x_i)\nonumber\\
&=&\big((x_i-a_i)(b_i-x_i)\widehat\phi_{i,1}(x_i),\cdots,(x_i-a_i)(b_i-x_i)
\widehat\phi_{i,p}(x_i)\big)^\top,
\end{eqnarray*}
where $\widehat\phi_i(x_i;\theta_i)$ is an FNN from $\mathbb R$ to $\mathbb R^p$ with sufficiently
smooth activation functions, such that $\Psi(x)\in H_0^1(\Omega)$.


\begin{theorem}\label{Error_estimate_Theorem}
If $\bar u$ and $u$ are solutions of (\ref{Laplace_Problem}) and (\ref{Laplace_Problem_Modified}), 
respectively. 
Assume the coefficients $a$ and $\widehat a$ has the same lower bound $a_{\min}$, 
both $b$ and $\widehat b$ are positive. 
Then 
\begin{eqnarray}\label{Error_Estimate}
\|\bar u-u\|_V \leq \left(\|f-\widehat f\|_{-1} +\left(\|a-\widehat a\|_{L^\infty(\Omega)}+
\|b-\widehat b\|_{L^\infty(\Omega)}\right)\|f\|_{-1}\right),
\end{eqnarray}
where $V:=H_0^1(\Omega)$ and the energy norm $\|\cdot\|_V$ is defined as follows 
\begin{eqnarray*}
\|v\|_V:=\sqrt{(\widehat a \nabla v, \nabla v)+(\widehat bv,v)},\ \ \ \forall v\in V. 
\end{eqnarray*}
\end{theorem}
\begin{proof}
The variational form for the equation (\ref{Laplace_Problem}) can be defined as follows
\begin{eqnarray}\label{Weak_Laplace_Problem}
(a\nabla\bar u,\nabla v)+(b\bar u,v)=(f,v),\ \ \ \forall v\in V.
\end{eqnarray}
Similarly, the weak form for the modified problem (\ref{Laplace_Problem_Modified}) 
can be given as follows
\begin{eqnarray}\label{Weak_Laplace_Problem_2}
(\widehat a\nabla u,\nabla v)+(\widehat b u,v)=(\widehat f,v),\ \ \ \forall v\in V.
\end{eqnarray}
From (\ref{Weak_Laplace_Problem}), (\ref{Weak_Laplace_Problem_2}) and the regularity 
$\|u\|_1 \leq C\|f\|_{-1}$, we have following 
estimates 
\begin{eqnarray}\label{Modified_Error_Estimate}
&&(\widehat a\nabla(\bar u- u),\nabla v)+(\widehat b (\bar u- u),v) 
=(f-\widehat f,v)+((\widehat a-a)\nabla u, \nabla v)+((\widehat b-b)u,v)\nonumber\\
&&\leq \left(\|f-\widehat f\|_{-1} +\|a-\widehat a\|_{L^\infty(\Omega)}\|u\|_1+
\|b-\widehat b\|_{L^\infty(\Omega)}\|u\|_0\right)\|v\|_1\nonumber\\
&&\leq C\left(\|f-\widehat f\|_{-1} +\left(\|a-\widehat a\|_{L^\infty(\Omega)}+
\|b-\widehat b\|_{L^\infty(\Omega)}\right)\|f\|_{-1}\right)\|v\|_1,\ \ \ \forall v\in V.
\end{eqnarray}
Setting $v=\bar u- u$ in (\ref{Modified_Error_Estimate}) leads to the following inequality 
\begin{eqnarray*}
\|\bar u - u\|_V \leq C \left(\|f-\widehat f\|_{-1} +\left(\|a-\widehat a\|_{L^\infty(\Omega)}+
\|b-\widehat b\|_{L^\infty(\Omega)}\right)\|f\|_{-1}\right).
\end{eqnarray*}
This is the desired result (\ref{Error_Estimate}) and the proof is complete. 
\end{proof}
Theorem \ref{Error_estimate_Theorem} shows that the way in this 
section is reasonable for solving partial differential equations 
with the coefficients being interpolated by the TNN functions. 
The detailed computing process can be decomposed into two steps. 
In the first step, we use TNN based machine learning method 
by Algorithm \ref{Algorithm_Interpolation} to interpolate the 
non-tensor-product-type coefficients or source term of the concerned 
problem (\ref{Laplace_Problem}). In the second step, 
the TNN based machine learning method \cite{WangLinLiaoLiuXie} 
is adopted to solve the modified problem (\ref{Laplace_Problem_Modified}).

\section{Numerical examples}\label{Section_Numerical_Examples}
In this section, we provide two numerical examples to validate the TNN interpolation method 
which is defined by Algorithm \ref{Algorithm_Interpolation}.  
In order to check the validation, the TNN interpolation method is used to 
compute the integrations of (possibly non-tensor-product-type) high dimensional functions and 
solve the high dimensional partial differential equations with 
non-tensor-product type of coefficients and source term.

\subsection{High dimensional integration}
In this subsection, we check the numerical performance of the TNN 
interpolation for the high dimensional integrations.  
For this aim, we set the desired function 
$f(x) = \exp\left(\sum_{i=1}^8 x_i^2\right)$ on the domain $\Omega=[0,1]^8$. 
The reason to choose a tensor-product-type function $f(x)$ is 
that we can compute the integration error accurately.

We execute the interpolation process following 
Algorithm \ref{Algorithm_Interpolation}. In each iteration, 
we uniformly sample $8000$ points in domain $[0,1]^8$. 
As for network structure, we choose the rank $p=50$, 
and each subnetwork of TNN is chosen as the FNN with two hidden 
layers and each hidden layer has 50 neurons. 
The sine function is chosen as the activation function. We conduct 
the optimization process using the Adam optimizer with learning rate 0.003 
in the first 50,000 steps and then the LBFGS optimizer with learning 
rate 0.1 in the subsequent 200 steps. We repeat the iterations 
for 20 times to obtain the final approximate function. In order to validate 
the accuracy of the integration using TNN interpolation method, 
the interval $[0,1]$ is decomposed into 100 subintervals and 16 
Gauss points is chosen on each subinterval for computing 
the integration error using the quadrature scheme in Section \ref{Section_Integration}. 
The error of the integration for $f(x)$ 
on $[0,1]^8$ using the appropriate function $\Psi_{\text{TNN}}(x)$ is as follows:
\[
\begin{aligned}
\left|\int_\Omega f(x)dx - \int_\Omega \Psi_{\text{TNN}}(x)dx\right| &\approx 8.813175\text{e{-07}}.
\end{aligned}    
\]
This error result shows that it is feasible to use the TNN interpolation  
to compute the high dimensional integrations. 
In the next subsection, we will combine the interpolation method 
with the TNN machine learning method to solve high dimensional elliptic equations.

\subsection{TNN interpolation for solving partial differential equation}
In this subsection, we consider the following Poisson boundary value problem: Find $\bar u\in H_0^1(\Omega)$ 
such that 
\vskip-0.2cm 
\begin{equation}\label{Example_2}
\left \{
\begin{aligned}
-\Delta \bar u &= f, \quad x \in \Omega, \\
\bar u &= 0, \quad x \in \partial \Omega,            
\end{aligned}
\right.
\end{equation}
where $\Omega = [-1,1]^d$, and the source term $f(x)$ is defined as follows  
\begin{equation*}
\begin{aligned}
f = \, & \exp \left(\prod_{i=1}^d (1+x_i)(1-x_i)\right) 
\cdot \sum_{k=1}^d\left(-4x_k^2 \prod_{i \neq k}(1+x_i)^2(1-x_i)^2 
+2\prod_{i \neq k}(1+x_i)(1-x_i) \right),
\end{aligned}
\end{equation*}
such that the exact solution is $\bar u(x) = \exp \left(\prod_{i=1}^d (1+x_i)(1-x_i)\right)-1$.  
It is easy to know that $g(\mathbf{x}) = \exp \left(\prod_{i=1}^d (1+x_i)(1-x_i)\right)$ 
is not of tensor product type in the source term $f(x)$.

In order to solve (\ref{Example_2}), we first do the TNN interpolation 
for the function $g(x)$. For this aim, we use a TNN function $g^{\text{TNN}}$ 
which has the following form
\begin{equation*}
g^{\text{TNN}} = \sum_{j=1}^{p}\alpha_j \prod_{i=1}^d \phi_{ij}(x_i),
\end{equation*}
to interpolate the function $g(x)$. 
We conduct the following iteration for 20 times. In each iteration, 
we uniformly sample $50,000$ points in domain $[-1,1]^d$. 
The TNN structure is built with the rank $p=100$ and each subnetwork being chosen 
as the FNN with two hidden layers and 50 neurons in each hidden layer. 
We use the sine function as the activation function. 
The Adam optimizer with learning rate 0.003 is adopted in the first 30,000 steps 
and then the LBFGS optimizer with learning rate 0.1 follows in the subsequent 2,000 steps. 
To validate the accuracy of the interpolation, we uniformly sample $N=10,000$ 
points $\{\mathbf{z}_k\}, k=1,\dots,10,000$ in $[-1,1]^d$ 
and compute the RMSE and $\ell_2$ relative error as follows 
\begin{equation*}
\text{RMSE:} \sqrt{\frac{\sum_{k=1}^N \left(g(\mathbf{z}_k) 
- g^{\text{TNN}}(\mathbf{z}_k)\right)^2}{N}},\ \   
\ell^2 \text{ relative error: } \sqrt{\frac{\sum_{k=1}^N \left(g(\mathbf{z}_k) - g^{\text{TNN}}(\mathbf{z}_k)\right)^2}{\sum_{k=1}^N \left(g(\mathbf{z}_k)\right)^2}}. 
\end{equation*}
on these test points.
Table \ref{Error_TNN_Interpolation_g} shows the corresponding errors 
of the TNN interpolation on the test points for $d=5,10,20$. 
\begin{table}[htb]
\begin{center}
\caption{The errors on the 10,000 points of the TNN interpolation for $d=5, 10, 20$}\label{Error_TNN_Interpolation_g}
\begin{tabular}{||c|c|c|c||}
\hline\hline
$d$ & 5 & 10 & 20 \\
\hline
RMSE & 1.4635e-06 & 8.5571e-07 & 2.1791e-07    \\
\hline
$\ell^2$ relative error & 1.2375e-06 & 8.3914e-07 & 2.1784e-07 \\ 
\hline\hline
\end{tabular}
\end{center}
\end{table}

Based on the TNN interpolation $g^{\text{TNN}}$, 
the equation (\ref{Laplace_Problem}) is modified by 
replacing $g(\mathbf{x})$ with $g^{\text{TNN}}$. 
Then we solve the following modified Poisson equation: 
Find $u\in H_0^1(\Omega)$ such that 
\begin{equation}\label{Equation_modified_poisson}
\left \{
\begin{aligned}
-\Delta u &= \widehat{f}, \quad x \in \Omega,\\
u &= 0, \quad x \in \partial \Omega,            
\end{aligned}
\right.
\end{equation}
where
\begin{equation*}
\widehat{f} = g^{\text{TNN}} \cdot 
\sum_{k=1}^d\left(-4x_k^2 \prod_{i \neq k}(1+x_i)^2(1-x_i)^2+2
\prod_{i \neq k}(1+x_i)(1-x_i) \right).
\end{equation*}
Then the TNN-based machine learning method (Algorithm 1 in \cite{WangLinLiaoLiuXie}) 
is adopted to solve (\ref{Equation_modified_poisson}) with the loss function being defined by
the a posteriori error estimators.  
We build a TNN with rank $p=100$ and each subnetwork is chosen as the FNN 
with three hidden layers and 100 neurons in each hidden layer. 
The Adam optimizer with learning rate 0.003 is adopted in the first 50,000 steps 
and then the LBFGS optimizer with learing rate 0.1 is used for the subsequent 10,000 steps. 
For the quadrature scheme, the interval is decomposed into 100 subintervals 
and 16 Gauss points are chosen on each subinterval. 
We compute the errors on test points to validate the accuracy. 
The corresponding errors between the TNN approximation $u^{\text{TNN}}$ 
and the exact solution $\bar u(x)$ on $N=50,000$ test points for $d=5,10,20$ 
are shown in Table \ref{Error_TNN_Solution}.
\begin{table}[htb]
\begin{center}
\caption{The errors on the 50,000 points of the TNN approximation for $d=5, 10, 20$}
\label{Error_TNN_Solution}
\begin{tabular}{||c|c|c|c||}
\hline\hline
$d$ & 5 & 10 & 20 \\
\hline
RMSE & 1.8232e-07 & 1.0972e-07 & 4.9389e-08    \\
\hline
$\ell^2$ relative error & 6.8637e-07 & 2.3012e-06 & 2.9410e-05 \\ 
\hline\hline
\end{tabular}
\end{center}
\end{table}

\section{Conclusion}
In this paper, we design a type of TNN based interpolation 
for high dimensional functions which has no tensor-product structure. 
Different from the general interpolation, the TNN interpolation 
has the tensor product structure and then the corresponding 
high-dimensional integration can be computed with high accuracy.
Then, benefit from the high accuracy of the high-dimensional integration for 
the TNN interpolation,  the TNN based machine learning method 
can be adopted to solve the high-dimensional differential equations, 
which has non-tensor-product-type coefficients and source term,  with high accuracy. 
We believe that the ability of TNN based interpolation and machine 
learning method will bring more applications in solving linear 
and nonlinear high-dimensional PDEs. These will be our future work.

\bibliographystyle{plain}
\bibliography{references}

\end{document}